\documentclass[11pt]{article}

\usepackage{amscd,amsmath,amssymb}

\def\eqref#1{(\ref{#1})}
\newcommand{\goth}{\frak}
\newcommand{\g}{{\frak g}}
\newcommand{\arrow}{{\:\longrightarrow\:}}
\newcommand{\Z}{{\Bbb Z}}
\newcommand{\C}{{\Bbb C}}
\newcommand{\R}{{\Bbb R}}

\newcommand{\restrict}[1]{{\left|_{{\phantom{|}\!\!}_{#1}}\right.}}

\renewcommand{\tilde}{\widetilde}
\renewcommand{\bar}{\overline}
\renewcommand{\phi}{\varphi}
\renewcommand{\epsilon}{\varepsilon}
\renewcommand{\geq}{\geqslant}
\renewcommand{\leq}{\leqslant}

\newcommand{\End}{\operatorname{End}}

\newcommand{\codim}{\operatorname{codim}}

\newcommand{\comment}[1]{{}}

\def\blacksquare{\hbox{\vrule width 4pt height 4pt depth 0pt}}
\def\endproof{\blacksquare}

\makeatletter

\@ifundefined{Bbb}
     {\newcommand{\Bbb}[1]{{\mathbb #1}}}%
{}

\newcommand{\ps@verbit}{%
  \renewcommand{\@oddhead}{%
          \scriptsize
          {Symplectic resolutions}
          \hfil\tiny {M. Verbitsky \ \ \ \ 30 Mar 1999}}
  \renewcommand{\@evenhead}{\@oddhead}
  \renewcommand{\@oddfoot}{\hfil\thepage\hfil}
  \renewcommand{\@evenfoot}{\@oddfoot}}
 
\newcounter{Mycounter}[section]
\newcounter{lemma}[section]
\renewcommand{\thelemma}{{Lemma \thesection.\arabic{lemma}}}
\newcommand{\lemma}{%
     \setcounter{lemma}{\value{Mycounter}}
     \refstepcounter{lemma}
     \stepcounter{Mycounter}
     {\bf \thelemma:\ }}

\newcounter{claim}[section]
\renewcommand{\theclaim}{{Claim \thesection.\arabic{claim}}}
\newcommand{\claim}{%
     \setcounter{claim}{\value{Mycounter}}
     \refstepcounter{claim}
     \stepcounter{Mycounter}
     {\bf \theclaim:\ }}

\newcounter{sublemma}[section]
\newcounter{corollary}[section]
\renewcommand{\thecorollary}{{Corollary \thesection.\arabic{corollary}}}
\newcommand{\corollary}{%
     \setcounter{corollary}{\value{Mycounter}}
     \refstepcounter{corollary}
     \stepcounter{Mycounter}
     {\bf \thecorollary:\ }}

\newcounter{theorem}[section]
\renewcommand{\thetheorem}{{Theorem \thesection.\arabic{theorem}}}
\newcommand{\theorem}{%
     \setcounter{theorem}{\value{Mycounter}}
     \refstepcounter{theorem}
     \stepcounter{Mycounter}
     {\bf \thetheorem:\ }}

\newcounter{conjecture}[section]
\renewcommand{\theconjecture}{{Conjecture \thesection.\arabic{conjecture}}}
\newcommand{\conjecture}{%
     \setcounter{conjecture}{\value{Mycounter}}
     \refstepcounter{conjecture}
     \stepcounter{Mycounter}
     {\bf \theconjecture:\ }}

\newcounter{proposition}[section]
\renewcommand{\theproposition}
       {{Proposition \thesection.\arabic{proposition}}}
\newcommand{\proposition}{%
     \setcounter{proposition}{\value{Mycounter}}
     \refstepcounter{proposition}
     \stepcounter{Mycounter}
     {\bf \theproposition:\ }}

\newcounter{definition}[section]
\renewcommand{\thedefinition}
       {{Definition \thesection.\arabic{definition}}}
\newcommand{\definition}{%
     \setcounter{definition}{\value{Mycounter}}
     \refstepcounter{definition}
     \stepcounter{Mycounter}
     {\bf \thedefinition:\ }}

\newcounter{example}[section]
\renewcommand{\theexample}{{Example \thesection.\arabic{example}}}
\newcommand{\example}{%
     \setcounter{example}{\value{Mycounter}}
     \refstepcounter{example}
     \stepcounter{Mycounter}
     {\bf \theexample:\ }}

\newcounter{remark}[section]
\renewcommand{\theremark}{{Remark \thesection.\arabic{remark}}}
\newcommand{\remark}{%
     \setcounter{remark}{\value{Mycounter}}
     \refstepcounter{remark}
     \stepcounter{Mycounter}
     {\bf \theremark:\ }}

\newcounter{problem}[section]
\newcounter{question}[section]
\renewcommand{\thequestion}{{Question \thesection.\arabic{question}}}
\newcommand{\question}{%
     \setcounter{question}{\value{Mycounter}}
     \refstepcounter{question}
     \stepcounter{Mycounter}
     {\bf \thequestion:\ }}

\@addtoreset{equation}{section}
\@addtoreset{footnote}{section}
\makeatother

\begin{document}

\begin{center}
{\Large\bf
Holomorphic symplectic geometry and orbifold singularities
}
\\[4mm]
Misha Verbitsky,\\[4mm]
{\tt verbit@@thelema.dnttm.rssi.ru}
\end{center}

{\small 
\hspace{0.2\linewidth}
\begin{minipage}[t]{0.7\linewidth}
Let $G$ be a finite group acting on a symplectic complex vector 
space $V$. Assume that the quotient $V/G$ has a holomorphic 
symplectic resolution. We prove that $G$ is generated by 
``symplectic reflections'', i.e. symplectomorphisms with 
fixed space of codimension 2 in $V$. Symplectic resolutions 
are always semismall. A crepant resolution of $V/G$ is 
always symplectic. We give a symplectic version of 
Nakamura conjectures.
\end{minipage}
}

\tableofcontents


\section{Introduction}


\subsection{Symplectic desingularizations in 
algebraic geometry and representation theory}
\label{_symple_in_alge_Subsection_}

Let $V$ be a complex vector space, 
and $G$ a finite group acting on $V$ 
by linear transformations. The variety $X=V/G$ is
usually singular, and this paper deals with its desingularizations
(also called resolutions). A resolution of
$X$ is a proper birational map
$\pi:\; \tilde X\arrow X$ such that
$\tilde X$ is smooth, and $\pi$ is an isomorphism
outside of singularities of $X$. 

\hfill

A singularity of the type $V/G$ is called {\bf a quotient}
or {\bf orbifold} singularity.

\hfill

The crepant resolutions of $X$ are resolutions
$\pi:\; \tilde X\arrow X$ such that the canonical class
of $\tilde X$ is obtained as a pullback of a canonical
class of $X$ (see \ref{_crepant_Definition_}). 

\hfill

\example\label{_Hilbert_Example_}
The Hilbert scheme of $n$ points on $\C^2$ provides a crepant
resolution of the quotient $(\C^2)^n/S_n$ of $(\C^2)^n$
by the natural action of the symmetric group $S_n$
(this is well known; see e. g. \cite{_Nakajima_Hilb_}).

\hfill

The crepant resolutions of quotient singularities in dimension 3 and more
became  a focus of intense study after the paper \cite{_Ito_Reid_}
of Y. Ito and M. Reid,
because of their relations with physics and with the 
theory of Hilbert schemes (\cite{_Ito_Nakajima_}).
For a history of these questions and their relevance to 
the mirror symmetry, see \cite{_Reid:McKay_} and \cite{_Batyrev_Dais_}.

\hfill

Another reason to study the crepant resolutions comes from
the holomorphic symplectic geometry and representation
theory. Suppose that
a complex vector space space $V$ is equipped with 
a $\C$-valued symplectic form, and $G$ acts on $V$
by symplectic transformations. The desingularization
$\tilde X\arrow X$ is called {\bf a symplectic resolution}
if $\tilde X$ is holomorphic symplectic, and the
holomorphic symplectic form on $\tilde X$ is lifted from
$X$ (see \ref{_symple_reso_Definition_}).
Clearly, the symplectic resolutions are always
crepant. It turns out that, conversely,
any crepant resolution of $(V/G)$ is 
symplectic (\ref{_crepant+>symplectic_Theorem_}).
Symplectic resolutions were studied 
by R. Bezrukavnikov and V. Ginzburg (1998, 
unpublished), who worked in the following situation. 
Consider a semi-simple Lie algebra $\g$ over $\C$, and
its Cartan subalgebra $\goth h \subset \g$.  
Bezrukavnikov and Ginzburg considered the space
${\goth h} \oplus {\goth h}^*$ 
(a direct sum of a Cartan
algebra with its dual). Clearly, ${\goth h} \oplus {\goth h}^*$ is 
a symplectic vector space, equipped with a natural action
of the Weyl group $W$ of $\g$. They suggested that
the variety ${\goth h} \oplus {\goth h}^*/W$ 
admits a natural symplectic desingularization,
and this desingularization is hyperk\"ahler. This is true
for the case $\g = {\goth{sl}}(n)$, because in this case
\[ ({\goth h} \oplus {\goth h}^*)\big /W \cong (\C^2)^n/S_n,
\]
and the desingularization is provided by the Hilbert
scheme (\ref{_Hilbert_Example_}).
The conjectural desingularizations of 
Bezrukavnikov and Ginzburg are quite important, because
they generalize the usual Hilbert schemes.

\hfill

A second example when Bezrukavnikov-Ginzburg conjecture
is valid was invented by A. Kuznetsov, who considered
the Lie algebras associated with the Dynkin diagrams
$B_n$, $C_n$. The corresponding symplectic desingularization
of ${\goth h} \oplus {\goth h}^*/W$ is birational to the Hilbert scheme of 
the total cotangent space $T^*\C P^1$ of $\C P^1$.
Kuznetsov's construction is explained in more details in
\cite{_KV:2_}. There is some indication that 
Bezrukavnikov-Ginzburg conjecture is not valid
for the Dynkin diagrams $G_2$, $D_n$ ($n\geq 4$)
and $E_n$ ($n = 6, 7, 8$).

\hfill

Another reason to study the symplectic desingularization comes 
from the hyperk\"ahler geometry. 
Consider a compact complex torus $T$,
$\dim_\C T=2$, and its $n$-th Hilbert scheme 
of points $T{[n]}$. Let $Alb:\; T^{[n]}\arrow T$ be the
Albanese map. A {\bf generalized Kummer variety} 
$K^{[n-1]}$  is defined as
\[ K^{[n-1]}\subset T^{[n]}, \ \ K^{[n-1]}:= Alb^{-1}(0).
\]
The variety $K^{[n-1]}$ is smooth and holomorphically symplectic 
(\cite{_Beauville_}). By Calabi-Yau theorem
(\cite{_Yau:Calabi-Yau_}, \cite{_Beauville_}), 
the variety $K^{[n]}$ is equipped with a set of hyperk\"ahler
structures, parametrized by the K\"ahler cone.
In \cite{_KV:1_}, it was falsely claimed that,
for a generic hyperk\"ahler structure, 
$K^{[n]}$ has no subvarieties compatible with the hyperk\"ahler
structure (such subvarieties are called trianalytic, 
see \cite{_Verbitsky:Symplectic_II_}). We mentioned above the
simple agrument used by Kuznetsov to prove the 
existence of symplectic desingularizations of 
${\goth h} \oplus {\goth h}^*/W$ for the Dynkin diagrams $B_n$, $C_n$ 
The same argument proves existence of trianalytic subvarieties
of generalized Kummer varieties (\cite{_KV:2_}, Theorem 6.10). 

In \cite{_KV:2_} (Section 4), 
this topic was pursued further. It turns out
that all trianalytic subvarieties of generalized
Kummer varieties (at least, for generic hyperk\"ahler 
structures) are isomorphic to symplectic desingularizations
of a quotient of a compact torus by an action of a
Weyl group. This establishes
a very interesting relation between 
the Dynkin diagrams and hyperk\"ahler 
geometry, and motivates the study
of symplectic desingularization of quotient
singularities.

\subsection[Symplectic desingularizations and symplectic reflections]
{Symplectic desingularizations \\ and symplectic reflections}

In this paper, we carry the argument used in
\cite{_KV:2_} a step further, to obtain information about 
the structure of finite groups $G\subset Sp(V)$ such that
$V/G$ admits a symplectic desingularization.
This is done as follows. Let $g\in \End(V)$ 
be a symplectomorphism of finite order. 
We say that $g$ is a {\bf symplectic reflection} if
\[ 
   \codim_V\bigg(\{x\in V\;\; |\; \; g(x) =x\}\bigg) =2,
\]
that is, the dimension of the fixed set of $g$ 
is maximal possible for non-trivial $g$ 
(see \ref{_symple_refle_Definition_}). 
This definition parallels that of complex reflections --
a complex reflection is an endomorphism of finite order
with fixed point set of codimension 1. The main result
of this paper is the following theorem.

\hfill

\theorem \label{_main_intro_Theorem_}
Let $V$ be a symplectic vector space over $\C$, and
$G\subset Sp(V)$ a finite group of symplecic transformations.
Assume that $V/G$ admits a symplectic resolution. Then
$G$ is generated by symplectic reflections.

\hfill

{\bf Proof:} This is \ref{_genera_symple_main_Theorem_}, 
which is proven in Section \ref{_main_proof_Section_}. 
\endproof

\hfill

This result is analogous to a well known theorem
(\ref{_comple_refle_Proposition_}, 
\ref{_quotie_by_comple_refle_smooth_Remark_}), stating
that for any finite group of endomorphisms
$G\subset \End V$, $V/G$ is smooth if and only
if $G$ is generated by complex reflections.
However, the ``if'' part in symplectic case
is not proven and is likely false
(see \ref{_Nakamura_Conjecture_}).

\hfill

In connection with \ref{_main_intro_Theorem_},
the following questions appear.

\hfill

\question
Is it possible to classify the groups 
generated by symplectic reflections?

\hfill

A complete classification of groups generated by complex reflections
was obtained in \cite{_Shepard_Todd_} 
(Shephard and Todd, 1954).

\hfill

The basic example of groups generated 
by symplectic reflections is the following.

\hfill

\example\label{_symple_refle_from_comple_refle_Example_}
Given an action of a group $G\subset \End W$ 
we consider the natural action of $G$ on $\End(W\oplus W^*)$.
Clearly, $G$ acts on $W\oplus W^*$ preserving 
the standard symplectic structure. The action
of $G$ on $W\oplus W^*$ is generated by symplectic reflections
if and only if the action of $G$ on $W$ is generated by complex
reflections (see the proof of \ref{_Kale_main_Theorem_}).
 
\hfill

The Weyl group acting on ${\goth h}\oplus{\goth h}^*$
(Subsection \ref{_symple_in_alge_Subsection_})
is a special case of \ref{_symple_refle_from_comple_refle_Example_}.

\hfill

\remark
Not all the groups generated by symplectic reflections
are provided by \ref{_symple_refle_from_comple_refle_Example_}.
Take, for instance, any finite subgroup $G\subset SL(V)$,
$\dim V =2$. Clearly, all non-trivial elements 
of $G$ are symplectic reflections. The group 
$G$ is obtained from \ref{_symple_refle_from_comple_refle_Example_}
if and only if $G$ preserves a direct decomposition
$V=W_1\oplus W_2$, $\dim W_1=\dim W_2=1$.
The finite subgroups of $SL(2)$ are well known,
and for most of them such a decomposition 
does not exist.

\hfill

Another question appearing in connection
with \ref{_main_intro_Theorem_}
is the following

\hfill

\question\label{_Ginzburg_Question_}
Let $G\subset \End V$ be a subgroup generated by symplectic
reflections. Determine whether $V/G$ admits a symplectic
resolution.

\hfill

In the case $\dim V=2$, the answer is ``always'' by the classical
results of Du Val; in the next non-trivial case
($\dim V =4$) the answer is unknown already.

\hfill

An example of a Weyl group of $G_2$ acting on 
${\goth h} \oplus {\goth h}^*\cong \C^4$ motivates
the following version of Nakamura's conjecture
(\cite{_Reid:McKay_}).

\hfill

\conjecture\label{_Nakamura_Conjecture_}
Let $G\subset Sp(V)$ be a finite group acting on a symplectic
$\C$-vector space, and $\tilde X:= \text{Hilb}^GV$ the
$G$-Hilbert scheme (\cite{_Reid:McKay_}). 
\begin{description}
\item[(i)] Then
$\tilde X$ is a smooth holomorphic symplectic 
orbifold with singularities in codimension $\geq 4$, and
the natural map $\pi:\; \tilde X \arrow V/G$ is an
orbifold desingularization of $V/G$. 
\item[(ii)] Moreover, for any crepant
orbifold desingularization $\pi_1:\; \tilde X_1\arrow V/G$
with $\dim Sing \tilde X_1\geq 4$, the orbifold
$\tilde X_1$ is diffeomorphic to $\tilde X$.
\end{description}

\hfill

The second part of \ref{_Nakamura_Conjecture_}
is motivated by a result of Huybrechts 
\cite{_Huybrechts_}: birational holomorphic 
symplectic compact manifolds are diffeomorphic.

\subsection{Contents}

\begin{itemize}

\item The present Introduction is independent from the rest of this paper.

\item In Section \ref{_symple_reso_Section_}, we define symplectic
desingularizations and state their main properties. Any crepant 
desingularization of a quotient singularity $V/G$,  $G\subset Sp(V)$
is symplectic (\ref{_crepant+>symplectic_Theorem_}). 
Simplectic resolutions are semismall (\ref{_semismall_Theorem_}).

\item In Section \ref{_symple_refle_Section_}, we state our main
result: if a symplectic desingularization of $V/G$ exists, then
$G$ is generated by symplectic reflections 
(\ref{_genera_symple_main_Theorem_}). We illustrate
this result by a proof of a classical theorem about 
groups generated by complex reflections.

\item In Section \ref{_main_proof_Section_},
we give the proof of our main theorem. The proof
is based on the following preliminary proposition
(\ref{_simply_conne_Theorem_}). Let $\tilde X\arrow X$ 
be a resolution of a quotient-type singularity $X=V/G$.
Then the manifold $\tilde X$ is simply connected.

\end{itemize}


\section{Symplectic resolutions}
\label{_symple_reso_Section_}


Let $X$ be an irreducible complex analytic variety. A proper 
morphism $\pi:\; \; \tilde X\arrow X$ is called {\bf a resolution 
of singularities}, or {\bf a desingularization} of $X$ 
if $\tilde X$ is smooth and connected, and $\pi$ is an 
isomorphism outside of the set of singular points of $X$.

\hfill

\definition \label{_crepant_Definition_}
In the above assumptions, let $X$ be a normal variety,
and $U$ a non-singular part of $X$. Since $\pi$ 
is an isomorphism over $U$, we may consider 
$U$ as a subset in $\tilde X$. Let $C\subset U$ be a divisor 
associated with the canonical class of $U$, and $\tilde C$ 
a closure of $C$ in $\tilde X$. We say that
$\pi$ is {\bf crepant} if the divisor $\tilde C$ 
lies in the canonical class of $\tilde X$. 

In other words, crepant resolutions are those which preserve
the canonical class.

\hfill

\definition
Let $V$ be a vector space, $G$ a finite group acting on
$V$ by linear transformations, and $B\subset V$ a $G$-invariant open ball.
In this case, we say that $B/G$ has 
{\bf orbifold singularities} or {\bf singularities
of quotient type}.

\hfill

\definition \label{_symple_reso_Definition_}
Let $\pi:\; \tilde X\arrow X$ be a resolution of a quotient-type
singularity $X= V/G$, where $V$ is a symplectic $\C$-vector space.
Assume that $G$ acts on $V$ by linear transformations preserving the
symplectic form: $G\subset Sp(V)$. 
Consider the natural symplectic form $\Omega_X$ on $X$,
defined outside of the singularities of $X$. Assume that
the pullback $\pi^* \Omega$  can be extended to a holomorphic
symplectic form on $\tilde X$. Then $\pi:\; \tilde X\arrow X$
is called {\bf a symplectic resolution} of $X$. 

\hfill

\remark
Clearly, symplectic resolutions are crepant.

\hfill

The following preliminary theorem establishes the relationship
between the symplectic and crepant resolutions. We shall not use it
in this paper, but prove it here to validate the notion of
a symplectic desingularization.

\hfill

\theorem\label{_crepant+>symplectic_Theorem_}
Let $X= V/G$ be a quotient of a symplectic 
vector space by a finite group $G\subset Sp(V)$. Assume that 
$\pi:\;\tilde X\arrow X$ is a crepant resolution
of singularities. Then $\pi$ is a symplectic resolution,
in the sense of \ref{_symple_reso_Definition_}.

\hfill

{\bf Proof:} Another proof of this theorem is given in
\cite{_Kaledin_}, (Proposition 3.2), and in \cite{_Beau:Symple_sing_}, 
(Proposition 2.4).

Let $U\subset X$ be the set of all points where 
$\pi:\; \tilde X \arrow X$ is smooth, and 
$\tilde U:= \pi^{-1}(U)$ the corresponding 
subset of $\tilde X$. Consider the symplectic form
$\Omega_X$ as a section of the sheaf of holomorphic
2-forms $\Omega^2_U$.
Let $\Omega_{\tilde U}:=\pi^* \Omega_X$ 
be its pullback. We need to show that 
the form $\Omega_{\tilde U}$ can be 
extended to the whole $\tilde X$.
Consider a smooth Hermitian metric on $\tilde X$, and let $h$ be the
corresponding metric on $\Omega^2 \tilde X$. Then
$\Omega_{\tilde U}$ is a holomorphic section of a Hermitian
vector bundle $\Omega^2_{\tilde U}$ over $\tilde U$. Such a section can
be extended to a section of a bundle over $\tilde X$ unless it has 
singularities on the complement $\tilde X\backslash \tilde U$. 
To prove that
$\Omega_{\tilde U}$ can be extended to $\tilde X$, it remains
to show that for any compact set $K\subset \tilde X$ 
for every $x\in \tilde U\cap K$, the Hermitian 
norm of $\Omega_{\tilde U}$ is bounded
\[ h_x(\Omega_{\tilde U}\restrict x) < C_K
\]
by some constant $C_K$ depending on $K$.

A symplectic form $\Theta$ in a Hermitian vector space $L$ can be 
naturally
represented in the form 
\[ \Theta = \sum_{i=0}^{n-1} \lambda_i z_{i+1}\wedge z_{i+2} \]
for some orthonormal basis $z_1, ..., z_{2n}$ in $L^*$, where
$\lambda_i$ are non-negative real numbers.
The numbers $\lambda_i$ are called {\bf the eigenvalues of $\Theta$};
the set of eigenvalues is defined canonically by the symplectic
form and the Hermitian form. Denote the eigenvalues by 
$\lambda_i(\Theta)$. Clearly, for any $x\in \tilde U$, we have

\[  h_x(\Omega_{\tilde U}\restrict x) 
    \leq C \cdot \max_i \lambda_i^2(\Omega_{\tilde U}\restrict x).
\]
for some constant $C$ depending only on $\dim X$.
Therefore, to show that $h(\Omega_{\tilde U})$ is bounded, it suffices
to show that the eigenvalues of $\Omega_{\tilde U}$ are bounded.

Since the manifold $U$ is holomorphic symplectic, 
its canonicall class is trivial. Denote by 
$\eta\in \Omega^{\dim X} U$ 
the trivializing section of the canonical class,
obtained as an $\frac{\dim X}{2}$-th power of $\Omega_X$. 
Since the resolution $\pi:\tilde X\arrow X$ is crepant,
the section $\eta$ can be extended to a global section
$\tilde\eta$ of the canonical class of $\Omega^{\dim X}\tilde X$.
Therefore, the Hermitian norm $h(\tilde\eta)$ is bounded over any
compact. On the other hand, for all $x\in \tilde U$, 
the norm $h(\tilde\eta)\restrict x$ 
is equal (up to a canonical constant)
to a product of all eigenvalues of 
$\Omega_{\tilde U}\restrict x$. Therefore,
$h(\tilde\eta)$ is bounded implies
$h(\Omega_{\tilde U})$ is bounded, unless
some eigenvalues of $\Omega_{\tilde U}\restrict x$
tend to zero as $x$ tends to the complement
$\tilde X\backslash \tilde U$. Therefore, to prove that
the  eigenvalues of $\Omega_{\tilde U}$ are bounded from 
above, it suffices to show that these eigenvalues are bounded from
below by some positive constant.

Consider the K\"ahler metric on $U\subset X$ obtained from
the flat metric on $V$. The map $\pi:\tilde X\arrow X$ is analytic,
and therefore, Lipschitz on compact subsets $K\subset \tilde X$
(the Lipschitz constant being given by the supremum of the
absolute value of the derivative). By definition of Lipschitz mappings,
for any compact set $K\subset \tilde X$, we may assume that,
after rescaling the metric, the map $\pi\restrict K$ is  
decreasing distances. 

The form $\Omega_U$ is parallel with respect to the natural flat
coordinates on $U$. Clearly, the eigenvalues of $\Omega_U$
are constant. Since $\pi:\; \tilde U\arrow U$
is decreasing distances, the eigenvalues of $\Omega_{\tilde U}=\pi^*\Omega_U$
are bounded from below by a positive constant.
As we have shown earlier, this implies that
these eigenvalues are bounded from above.  
\ref{_crepant+>symplectic_Theorem_} is proven.
\endproof

\hfill

\hfill

\definition
Let $\pi:\; \tilde X\arrow X$ be a resolution of singularities.
The map $\pi$ is called {\bf semismall} if $X$ admits a stratification
$\goth S$ with open strata $U_i$, such that
\[ 
   \forall x\in U_i \;\; |\;\; \dim \pi^{-1}(x)\leq \frac{1}{2}\codim U_i
\]

\hfill

\definition \label{_G_strati_Definition_}
Let $V$ be a symplectic vector space, $G\subset Sp(V)$ a finite
group, $X= V/G$. For every subgroup $G_1\subset G$, consider its 
fixed point set $V_{G_1}$. Let $X_{G_1}\subset X$ 
be the image of $V_{G_1}$ under the quotient map 
$\sigma:\; V\arrow V/G$. Consider a stratification of 
$X$ with closed strata $X_{G_1}$, numbered by the 
subgroups $G_1\subset G$. This stratification is called
{\bf the $G$-stratification} of $X$. There is a
similar stratification of $V$, which is also called
a $G$-stratification.

\hfill

The main tool of our arguments is the following theorem.

\hfill

\theorem \label{_semismall_Theorem_}
Let $\pi:\; \tilde X\arrow X$ be a symplectic resolution of
a quotient singularity $X=V/G$, $G\in Sp(V)$. 
Then $\pi$ is semismall with respect to the $G$-stratification.

\hfill

{\bf Proof:} This statement easily follows from 
Proposition 4.16 and Proposition 4.5 of \cite{_Verbitsky:HilbertK3_}
(see also \cite{_Kaledin_}, Proposition 4.4).
\endproof


\section{Symplectic and complex reflections}
\label{_symple_refle_Section_}


\subsection{The statement of the main result}

Let $V$ be a symplectic $\C$-vector space.

\hfill

\definition\label{_symple_refle_Definition_}
Let $\gamma\in Sp(V)$ be an endomorphism of finite order.
We say that $\gamma$ is a {\bf symplectic reflection} if
$\codim_V V_\gamma=2$, where $V_\gamma$ is the
space of all vectors fixed by $\gamma$.

\hfill

This is a ``symplectic analogue'' of the usual notion of a complex
reflection. The complex reflection is a linear automorphism of 
a vector space fixing a subspace of codimension 1. Since the fixed
space of a symplectomorphism must be symplectic, the minimal
codimension of $V_\gamma$ is 2; the endomorphism
$\gamma$ is a symplectic reflection
when this minimum is reached.

We say that a group $G\subset Sp(V)$ is {\bf generated by 
symplectic reflection} if there exist
symplectic reflections $\gamma_1, ... \gamma_n\in G$
which generate $G$. 

\hfill

The main result of this paper is the following

\hfill

\theorem \label{_genera_symple_main_Theorem_}
Let $V$ be a symplectic vector space
over $\C$, and $G\subset Sp(V)$ a finite group. Assume that
the quotient $X=V/G$ admits a symplectic resolution. Then $G$
is generated by symplectic reflections. 

\hfill

The rest of this paper is taken by the proof of
\ref{_genera_symple_main_Theorem_}.

\hfill

As a first corollary of \ref{_genera_symple_main_Theorem_},
we obtain a new proof of the following theorem of Kaledin
(\cite{_Kaledin_}, Theorem 1.7).

\hfill

\theorem \label{_Kale_main_Theorem_}
Let $W$ be a complex vector space, $G\subset \End W$ a finite
group of endomorphisms of $W$, and $V:= W\oplus W^*$ 
be a symplectic space, obtained as a direct sum
of $W$ and its dual. Consider the natural embedding
$G\hookrightarrow Sp(V)$. Assume the $V/G$ has a symplectic
resolution. Then the action of $G$ on $W$ is generated by
complex reflections.

\hfill

{\bf Proof:} By \ref{_genera_symple_main_Theorem_},
the action of $G$ on $V$ is generated by symplectic reflections.
Take a symplectic reflection $g\in G$, and let $W_g$, $V_g$
be the fixed subspaces of the action of $g$ on $V_g$, $W_g$.
By definition, $V_g = W_g \oplus W_g^*$. Since $\codim V_g =2$,
we have $\codim W_g=1$. Therefore, $g$ acts on $W$
as a complex reflection. \endproof

\hfill

Using the arguments of 
\cite{_KV:2_}, Theorem 5.6, 
one can immediately obtain the following
corollary of \ref{_Kale_main_Theorem_}, which generalizes
\cite{_KV:2_}, Theorem 5.6.

\hfill

\corollary
Let $T$ be a 2-dimensional compact complex torus,
which is Mumford-Tate generic (see Definition 5.4 of 
\cite{_KV:2_}). Consider the natural holomorphic 
symplectic structure on $T$ and its $n$-th power $T^n$.
Let $G$ be a finite group acting on $T^n$ by 
symplectomorphisms and fixing a point 
\[ \hat x\in T^n, \ \ 
   \hat x=(\underbrace{x, ..., x}_{\text{$N$ times}}), 
\ \ x\in T.
\]
Assume that $T^n/G$ admits a symplectic resolution.
Then $G$ is a Weyl group associated with some reductive
Lie group $\g$. Moreover, the tangent space $T^n_{\hat x}$ is identified
as a representation of $G$ with ${\goth h} \oplus {\goth h}$,
where $\goth h$ is the Cartan algebra of $\g$.

\endproof

\subsection{Groups generated by complex reflections}

The proof of \ref{_genera_symple_main_Theorem_} 
is based on the same ideas as the proof
of the following well-known statement 
(\cite{_Bourbaki:Lie4-6_}, Ch. V, \S 5 Theorem 4).

\hfill

\proposition\label{_comple_refle_Proposition_}
Let $V$ be a complex vector space, and 
$G\subset GL(V)$ a finite group acting on $V$.
Assume that $X:= V/G$ is smooth. Then $G$ is generated by
complex reflections. 

\hfill

{\bf Proof:} Let $G_0\subset G$ be a maximal subgroup of $G$ 
generated by complex reflections. Clearly, $G_0$ is a normal 
subrgroup of $G$. Consider the corresponding quotient
 $X_0:= V/G_0$, and let
\[ \begin{CD}\tau:\;\; X_0 @>{*\big/(G/G_0)}>> X
\end{CD}
\]
be the natural quotient map. 

Since $G_0$ contains all the complex reflections, for
all $\eta$ in the complement $G\backslash G_0$, the fixed point
set $V_\eta\subset V$ has codimension $>1$.
Let
\[ S:= \bigcup\limits_{\eta\in G\backslash G_0} V_\eta, 
\]
and $S_0:= S/G_0\subset X_0$ its image in $X_0$. 
Clearly,
\begin{equation}\label{_codim_S_0_Equation_}
\codim_{X_{0}} S_0 >1.
\end{equation}

\hfill

The following claim is trivial.

\hfill

\claim \label{_free_act_Claim_} 
The group $\Gamma:= G/G_0$ acts freely on $X_0\backslash S_0$

\hfill

{\bf Proof:} Let $\gamma\in \Gamma$ be a non-trivial element, 
and $\tilde \gamma$ an element of $G\backslash G_0$ which is
mapped to $\eta$ under the natural quotient map. Consider a 
fixed point $x\in V/G_0$ of $\gamma$, and let $v\in V$ be a
point mapped to $x$ under the natural quotient map. Then,
$\tilde \gamma(v) = g(v)$, where $g\in G_0$. We obtain that
$v$ is a fixed point of $g^{-1}\tilde\gamma$.
Since $g^{-1}\tilde\gamma$ belongs to $G\backslash G_0$,
we have
\[ v \in  V_{g^{-1}\tilde\gamma}\subset 
   S= \bigcup\limits_{\eta\in G\backslash G_0} V_\eta.
\]
Therefore, $\gamma$ has no fixed point outside of $S_0=S/G_0$.
\endproof

\hfill

\remark\label{_free_act_Remark_}
Notice that in the proof of \ref{_free_act_Claim_} we
nowhere used the exact nature of the group $G_0$.
This means that \ref{_free_act_Claim_} holds for any
normal subgroup $G_0\subset G$: the quotient
$\Gamma:= G/G_0$ acts freely on the set 
$X_0\backslash S_0$ defined as above.

\hfill

Let $S_1:= S/G\subset X$. By \ref{_free_act_Claim_}, 
the natural quotient map
\[ X_0 \arrow X
\]
is etale over $X\backslash S_1$. By \eqref{_codim_S_0_Equation_},
$\codim S_1>1$. Since $X$ is smooth and
$\codim S_1>1$, the open embedding
$X\backslash S_1\hookrightarrow X$ induces an isomorphism on
the fundamental group:
\begin{equation}\label{_pi_1(X)_same_Equation_}
\pi_1(X) = \pi_1(X\backslash S_1)
\end{equation}

\hfill

Consider the scaling map $h_\lambda:\; V \arrow V$, $\lambda\in \C$:
\[ h_\lambda(v) = \lambda\cdot v.
\]
Clearly, this map is $G$-equivariant, and hence, 
can be extended to $X = V/G$. 
Taking all $\lambda\in \R$, we obtain a map
\[ h:\; X\times \R \arrow X, \ \ \ \ h(x, t) = h_t(x).
\]
This map establishes a contraction of $X$ into a point.
Therefore, $X$ is contractible, and $\pi_1(X)=0$.
By \eqref{_pi_1(X)_same_Equation_}, $\pi_1(X\backslash S_1)=0$.
This implies that obtained above etale covering
\[ \begin{CD}\tau:\;\; X_0 \backslash S_0@>{*\big/\Gamma}>> X\backslash S_1
\end{CD}
\]
is trivial. This covering is by construction
a Galois covering with the Galois group $\Gamma$;
hence, $\Gamma$ is trivial. This proves
\ref{_comple_refle_Proposition_}.
\endproof

\hfill

\remark \label{_quotie_by_comple_refle_smooth_Remark_}
A converse statement is also true
(\cite{_Bourbaki:Lie4-6_}, Ch. V, \S 5 Theorem 4).
Namely, let $X$ be a quotient of a vector space by an action of a group
generated by complex reflections. Then $X$ is smooth.


\section[Symplectic desingularizations and symplectic reflections]
{Symplectic desingularizations\\ and symplectic reflections}
\label{_main_proof_Section_}


In this section, we prove \ref{_genera_symple_main_Theorem_}.

\subsection{Fundamental groups of resolutions
of quotient singularities}

We use the following theorem, which seems to be well known.

\hfill

\theorem \label{_simply_conne_Theorem_}
Let $V$ be a linear space, and $G\subset GL(V)$ a finite group of
linear transformations. Consider the quotient $X:= V/G$, and let
$\pi:\; \tilde X\arrow X$ be a resolution of singularities. Then
the manifold $\tilde X$ is simply connected.

\hfill

{\bf Proof:} The following proof was suggested by F. Bogomolov.

\hfill

First of all, we construct a canonical surjection
\begin{equation}\label{_surje_on_tilde_X_Equation_}
 G \stackrel p \arrow \pi_1(\tilde X).
\end{equation}
This is done as follows. 
Let $Z\subset V$ be the union of all $G$-strata of codimension
$\geq 2$ (\ref{_G_strati_Definition_}).
For any $x\in V\backslash Z$, the stabilizing subgroup
$St_x(G)$ is generated by complex reflections. 
A quotient of $\C^n$ by an action of a 
group generated by complex reflection is smooth
(\ref{_quotie_by_comple_refle_smooth_Remark_}). 
Therefore, the quotient $(V\backslash Z)/G$ is smooth.

Since the map $\pi:\; \tilde X\arrow X$ is idenity
over smooth points of $X=V/G$, we have a canonical
open embedding
\[ (V\backslash Z)/G \hookrightarrow \tilde X.
\]
The following claim is trivial.

\hfill

\claim
Let $U\hookrightarrow Y$ be a Zariski open subvariety
of a complex variety $Y$. Then we have a natural 
epimorphism
\begin{equation}\label{_open_surje_on_pi_1_Equation_}
\pi_1(U) \arrow \pi_1(Y).
\end{equation}
\endproof

\hfill

This gives a natural surjective homomorphism
\begin{equation}\label{_pi_1_V_/G_to_tilde_X_Equation_}
\pi_1\bigg((V\backslash Z)/G\bigg) \arrow \pi_1(\tilde X).
\end{equation}

\hfill

Since $\codim_V Z>1$, we have $\pi_1(V\backslash Z) = \pi_1(V)$.
Therefore, the manifold $V\backslash Z$ is simply connected.

\hfill

\claim\label{_simply_conne_quotie_Claim_}
A quotient of a simply connected manifold $Y$ by an action
of a finite group $G$ has a fundamental group which is
a quotient of $G$.

\hfill

{\bf Proof:}
Any covering $\tilde Y \stackrel a\arrow Y/G$ 
correspond uniquely to a map $Y\stackrel b\arrow \tilde Y$,
in such a way that 
\[ a\circ b:\; Y\arrow Y/G\]
 is the quotient map.
\endproof

\hfill

By \ref{_simply_conne_quotie_Claim_}, 
the group $\pi_1\bigg((V\backslash Z)/G\bigg)$
is a quotient group of $G$. Using 
\eqref{_pi_1_V_/G_to_tilde_X_Equation_},
we obtain the surjection
\eqref{_surje_on_tilde_X_Equation_}.

\hfill

For any $g\in G$, denote the corresponding cyclic group
by $\langle g\rangle$. We extend the action of
$\langle g\rangle\subset \C^*$ to the action of $\C^*$
by the standard linear-algebraic argument.
Consider a general $\C^*$-orbit
$C_g\subset V$. Then
\begin{description}
\item[(i)]  $\langle g\rangle$ fixes $C_g$, and we have 
an embedding $\langle g\rangle\hookrightarrow \End(C_g)$
\item[(ii)]  $\langle g\rangle$ is the group $St(C_g)$
of all $\gamma\in G$
which map the line $C_g$ to itself.
\end{description}
Clearly, $C_g$ lies in $(V\backslash Z)/G$.
As a complex manifold, $C_g$ is naturally
isomorphic to $\C^*$. Denote by $\gamma_g$ the generator
of $\pi_1(C_g)\cong \Z$.

\hfill

\lemma\label{_C_g_gene_p_1_Lemma_}
Consider the natural group homomorphism
$\pi_1(C_g) \stackrel\tau\arrow \pi_1(V\backslash Z)/G)$.
Then $\tau$ maps $\gamma_g$ to the element $p(g)$ 
corresponding to $g\in G$ under the epimorphism 
\eqref{_surje_on_tilde_X_Equation_}.

\hfill

{\bf Proof:} Clear from 
\eqref{_surje_on_tilde_X_Equation_}.
\endproof

\hfill

Return to the proof of \ref{_simply_conne_Theorem_}.
Let $g\in G$. Consider the corresponding curve 
$C_g\stackrel i \arrow \tilde X$.
The map $\tilde X\arrow X=V/G$ is proper. Using valuative criterion
of properness, we extend the embedding $C_g\stackrel i \arrow \tilde X$
to a map $\bar C_g\stackrel i \arrow \tilde X$, where $\bar C_g \cong \C$
is the completion of $C_g$ in zero. Since $\bar C_g$ is simply
connected, the path $\gamma_g\subset C_g$ is contractible
in $\bar C_g\subset \tilde X$.
This proves \ref{_simply_conne_Theorem_}.
\endproof

\subsection{A subgroup generated by symplectic reflections}

Return to the assumptions of \ref{_genera_symple_main_Theorem_}.
Let $G_0\subset G$ be a subgroup of $G$ 
generated by all symplectic reflections in $G\subset \End(V)$.
Clearly, $G_0$ is a normal subgroup of $G$. 
We denote by $\tilde X_0$ the fibered product
\[ \tilde X_0:= (V/G_0)\times_{V/G} \tilde X,
\]
associated with the Cartesian square
\begin{equation}\label{_commu_square_Equation_}
\begin{CD}
\tilde X_0 @>{\tilde \sigma_0}>> \tilde X\\
@V{\pi_0}VV @VV\pi V \\
V/G_0 @>{\sigma_0}>> V/G. 
\end{CD}
\end{equation}
The top horizontal arrow 
$\tilde X_0 
\stackrel {\tilde \sigma_0}\arrow \tilde X$ 
is finite, since the bottom horizontal arrow is finite.
Denote the ramification variety of the map 
${\tilde \sigma_0}$ by $R\subset \tilde X$.

\hfill

\lemma \label{_codim_R_Lemma_}
We have $\codim_{\tilde X} R\geq 2$

\hfill

{\bf Proof:} Let $R_X\subset X$ be the ramification
variety for $\sigma_0:\;\; V/G_0 \arrow V/G$.
Then $R\subset \pi^{-1}(R_V)$.
Clearly, $R_X$ is a union of several strata of 
the $G$-stratification (\ref{_G_strati_Definition_}).
By \ref{_semismall_Theorem_},
for any stratum $U_i$ of a $G$-stratification,
we have
\[ 
   \codim \pi^{-1}(U_i) \geq \frac{1}{2} \codim U_1
\]
Therefore, to prove $\codim_{\tilde X} R\geq 2$,
it suffices to show that
\begin{equation}\label{_codim_R_V_Equation_}
\codim_X R_X \geq 4.
\end{equation}

\hfill

Let $Z\subset V$ be the union of fixed sets
of all $g\in G\backslash G_0$. By \ref{_free_act_Remark_}, the group
$(G/G_0)$ acts on $(V\backslash Z)/G_0$ without fixed points. Therefore,
the natural quotient map
\[ (V\backslash Z)/G_0 \stackrel{\sigma_0}\arrow (V\backslash Z)/G \]
is etale. By definition, $R_X$ is the ramification
variety for the map $\sigma_0:\;\; V/G_0 \arrow V/G$.
Therefore, $R_X\subset \sigma(Z)$,
where  $\sigma:\; V \arrow X$ is the quotient map.
Since the map $\sigma$ is finite, to prove 
\eqref{_codim_R_V_Equation_} it remains to show
that $\codim Z\geq 4$.
By definition,
\[ Z= \bigcup\limits_{\eta\in G\backslash G_0} V_\eta, 
\]
where $V_\eta$ denotes a fixed space of $\eta\in G\backslash G_0$.
To prove $\codim Z\geq 4$, we need to show 
that $\codim V_\eta\geq 4$.
The space $V_\eta$ is symplectic, hence its dimension
(and codimension) is even. On the other hand,
$\codim V_\eta>2$, because $\eta$ is {\it not}
a symplectic reflection. This proves 
\eqref{_codim_R_V_Equation_}. 
We proved \ref{_codim_R_Lemma_}.
\endproof

\hfill

\corollary
The natural embedding $\tilde X\backslash R\hookrightarrow \pi_1(\tilde X)$
induces an isomorphism
\[ \pi_1(\tilde X\backslash R) \cong \pi_1(\tilde X) \]

{\bf Proof:} The manifold $\tilde X$ is smooth, and
$\codim _{\tilde X}R >1$ by  \ref{_codim_R_Lemma_}.
\endproof

\hfill

By \ref{_simply_conne_Theorem_}, the manifold $\tilde X$
is simply connected. We obtain that 
$\tilde X\backslash R$ is also simply connected.
On the other hand, the map $\tilde \sigma_0$ of
\eqref{_commu_square_Equation_} induces a Galois
covering with the Galois group $G/ G_0$:
\[ \tilde \sigma_0:\; 
   \tilde \sigma_0^{-1}(\tilde X\backslash R)\arrow
    \tilde X\backslash R.
\]
Since $\tilde X\backslash R$ is simply connected, this
map is trivial, and its Galois group $G/ G_0$
is also trivial. This implies that $G$ coincides
with its subgroup $G_0$ generated by the symplectic
reflections. We obtained that $G$ is generated by symplectic
reflections. \ref{_genera_symple_main_Theorem_} is proven.
\endproof

\hfill

\hfill

{\bf Acknowledgements:}
This paper appeared in result of 
fruitful collaborative work with D. Kaledin.
F. Bogomolov was very helpful answering my questions; 
he also suggested the proof of simple-connectedness of $\tilde X$
(\ref{_simply_conne_Theorem_}). 
I am grateful to V. Batyrev, M. Leenson, A. Kuznetsov, M. 
Reid and T. Pantev for insightful comments.

\hfill


\end{document}